\documentclass[11pt]{amsart}

\theoremstyle{plain}
\newtheorem{theorem}                 {Theorem}      [section]
\newtheorem{proposition}  [theorem]  {Proposition}
\newtheorem{corollary}    [theorem]  {Corollary}
\newtheorem{lemma}        [theorem]  {Lemma}

\theoremstyle{definition}
\newtheorem{example}      [theorem]  {Example}
\newtheorem{remark}       [theorem]  {Remark}
\newtheorem{definition}   [theorem]  {Definition}

\numberwithin{equation}{section}

\def \1{\mbox{${\mathbf 1}$}}
\def \r{\mbox{${\mathbb R}$}}
\def \s{\mbox{${\mathbb S}$}}

\def \f{\mbox{$\phi$}}

\DeclareMathOperator{\trace}{trace}
\DeclareMathOperator{\grad}{grad}

\DeclareMathOperator{\ricci}{Ricci}

\begin{document}
\title{Biharmonic properties and conformal changes}
\author{A.~Balmu\c s}
\address{Faculty of Mathematics\\"Al.I. Cuza"
University of Ia\c si\\ Bd. Carol I, no. 11\\700506 Ia\c
si\\ROMANIA} \email{obalmus@uaic.ro, balmus1@yahoo.com}


\thanks{Partially supported by the Grant A.T. 73/2004, C.N.C.S.I.S, Romania}
\begin{abstract}
We construct a new class of biharmonic maps, which are the
critical points for the bienergy functional, by deforming
conformally the codomain metric of harmonic Riemannian submersions
such that they become nonharmonic but biharmonic.
\end{abstract}
\maketitle

\section{Introduction}

Let $\f:(M,g) \to (N,h)$ be a smooth map between Riemannian
manifolds. Its {\it tension field} is given by
$\tau(\phi)=\trace\nabla d\phi$ and, for a compact domain
$K\subseteq M$, the {\it bienergy} of $\f$ is
$$
E_2(\phi)=\frac{1}{2}\int_K |\tau(\phi)|^2v_g.
$$
The critical points of $E_2$ are called {\it biharmonic} maps.
 G.Y.~Jiang obtained in
\cite{GYJ1}, \cite{GYJ2} the first and the second variation
formula for the bienergy, showing that $\phi$ is biharmonic if and
only if
$$\tau_2(\phi)=-J(\tau(\phi))=0,$$
where $J=\triangle^\phi+\trace R^N(d\phi,\cdot)d\phi$ is the
Jacobi operator of $\phi$,
$\triangle^\phi=-\trace_g(\nabla^{\phi})^2$ is the Laplacian on
the sections of the pull-back bundle $\f^{-1}TN$ and $R^N$ is the
curvature operator
$R^N(X,Y)=[\nabla^N_X,\nabla^N_Y]-\nabla^N_{[X,Y]},\, \forall
X,Y\in C(TN).$

Since the Jacobi operator is linear, harmonic maps are biharmonic,
so we are interested in nonharmonic biharmonic maps.

There are results concerning the non-existence of nonharmonic
biharmonic maps when the sectional, or Ricci curvature of $N$ is
nonpositive \cite{GYJ1}, \cite{GYJ2}, \cite{CO1}. B.Y.~Chen and
S.~Ishikawa proved in \cite{CY} that there are no nonminimal
biharmonic submanifolds of the Euclidian $3$-dimensional space
$\r^3$ and, similarly, in \cite{RCSMCO2}, R.~Caddeo, S.~Montaldo
and C.~Oniciuc proved that there are no such submanifolds in
3-dimensional manifolds of negative constant sectional curvature
$-1$.

The nongeodesic biharmonic curves of surfaces of revolution in
$\r^3$ were studied in \cite{RCSMPP}, \cite{RCSMPP1}. In
\cite{RCSMCO1} the authors classified the nonminimal biharmonic
submanifolds of the unit Euclidian 3-dimensional sphere $\s^3$,
manifold of positive constant sectional curvature, and in
\cite{RCSMCO2} they gave two methods for construction of
nonminimal biharmonic submanifolds of $\s^n, n>3$.

A natural way of constructing nonharmonic biharmonic maps is the
following

\noindent {\it Given $\f:(M,g)\to(N,h)$ a harmonic map, one can
conformally change the metric $g$, or $h$, in order to get $\f$
nonharmonic but biharmonic.}

The first problem, i.e. the conformal change of the domain metric,
was studied by P.~Baird and D.~Kamissoko in \cite{PBDK}. They
called a metric $\widetilde{g}$, which renders the identity map
$\widetilde{\1}:(M,\widetilde{g})\to(M,g)$ biharmonic, a {\it
biharmonic metric} and proved that on Einstein manifolds the
conformally equivalent biharmonic metrics are defined only by
isoparametric functions.

In \cite{CO2}, \cite{CO3}, the author studied a particular case of
the second problem, i.e. the conformal change of the codomain
metric, providing a new method of constructing nonminimal
biharmonic submanifolds in spheres. The equatorial hypersphere
$\s^{n-1}$ is totally geodesic, thus minimal, in $\s^n$ and under
a certain conformal change of the Euclidean metric of $\s^n$, it
becomes nonminimal but biharmonic.

The goal of this paper is the study of the second problem, but,
instead of considering harmonic Riemannian immersions like in
\cite{CO2}, we treat the case of Riemannian submersions with
minimal fibres, or, equivalently, of harmonic Riemannian
submersions (see \cite{ER}). The first result is the reduction of
the second problem for harmonic Riemannian submersions to the
study of conformal changes with the property that the identity map
$\1:(N,h)\to(N, e^{2\rho}h)$ is nonharmonic and biharmonic. The
results concerning the biharmonicity of the identity after a
conformal change of the codomain  metric in the case of Einstein
spaces are similar to those obtained in \cite{PBDK}, but we
underline the fact that, in general, $\1:(N,h)\to(N, e^{2\rho}h)$
biharmonic is not equivalent to
$\widetilde{\1}:(N,e^{2\rho}h)\to(N,h)$ biharmonic.

\subsection*{Acknowledgement} The author warmly thanks Prof.
V.~Oproiu and C.~Oniciuc for useful discussions and ideas.

\section{Harmonic Riemannian submersions and conformal changes}

We have the following result

\begin{proposition}
Let $\pi:(M,g)\to (N,h)$ be a harmonic Riemannian submersion and
$\widetilde{h}$ a Riemannian metric on $N$. Denote by
$\widetilde{\pi}=\1\circ\pi:(M,g)\to(N,\widetilde{h})$, where
$\1:(N,h)\to (N,\widetilde{h})$ is the identity map. Then
$$\tau(\widetilde{\pi})=\tau(\1) \qquad \textrm{and}\qquad
\tau_2(\widetilde{\pi})=\tau_2(\1).$$
\end{proposition}

\begin{proof}
Let $\{X_a\}_{a=1}^n$ be a local orthonormal frame field on $N$,
then $\overline{X}_a=(X_a)^H$ is a local orthonormal base in
$T^HM$. We can complete $\{\overline{X}_a\}_{a=1}^n$ with
$\{\overline{X}_\alpha\}_{\alpha=n+1}^m\subset T^VM$  to a local
orthonormal frame field on $M$.

Using the chain rule for the tension field we obtain
\begin{eqnarray*}
\tau(\widetilde{\pi})&=& d\1(\tau(\pi))+\trace_g\nabla
d\1(d\pi\cdot,d\pi\cdot)=\trace_g\nabla d\1(d\pi\cdot,d\pi\cdot)
\\&=&\sum_{k=1}^m\nabla d\1
(d\pi(\overline{X}_k),d\pi(\overline{X}_k))=\sum_{a=1}^n\nabla
d\1(X_a,X_a) =\trace_h\nabla d\1\\&=& \tau(\1).
\end{eqnarray*}

In order to study the biharmonicity of $\widetilde\pi$ we compute
\begin{eqnarray*}
-\Delta\tau(\widetilde\pi)&=& \trace_g\nabla^{\widetilde{\pi}}
\nabla^{\widetilde{\pi}}\tau(\widetilde\pi)
\\&=&\sum_{k=1}^m\big\{\nabla^{\widetilde{\pi}}_{\overline{X}_k}
\nabla^{\widetilde{\pi}}_{\overline{X}_k}\tau(\widetilde\pi)
-\nabla^{\widetilde{\pi}}_{\nabla^M_{\overline{X}_k}
\overline{X}_k}\tau(\widetilde\pi)\big\}\\&=&\sum_{a=1}^n\big\{\widetilde{\nabla}_{X_a}
\widetilde{\nabla}_{X_a}\tau(\widetilde\pi)
-\widetilde{\nabla}_{\nabla^N_{X_a}
X_a}\tau(\widetilde\pi)\big\}-\sum_{\alpha=n+1}^m
\widetilde{\nabla}_{d\pi(\nabla^M_{\overline{X}_\alpha}
\overline{X}_\alpha)}\tau(\widetilde\pi)
\end{eqnarray*}

\begin{eqnarray*}
&=&\sum_{a=1}^n\big\{\widetilde{\nabla}_{X_a}
\widetilde{\nabla}_{X_a}\tau(\widetilde\pi)
-\widetilde{\nabla}_{\nabla^N_{X_a} X_a}\tau(\widetilde\pi)\big\}+
\widetilde{\nabla}_{d\1(\tau(\pi))}\tau(\widetilde\pi)\\
&=&\sum_{a=1}^n\big\{\widetilde{\nabla}_{X_a}
\widetilde{\nabla}_{X_a}\tau(\1)
-\widetilde{\nabla}_{\nabla^N_{X_a} X_a}\tau(\1)\big\}
\\&=&-\Delta\tau(\1)
\end{eqnarray*}
and
\begin{eqnarray*}
\trace_g\widetilde{R}(d\widetilde{\pi}\cdot,
\tau(\widetilde{\pi}))d\widetilde{\pi}\cdot&=&
\sum_{a=1}^n\widetilde{R}(X_a,\tau(\widetilde{\pi}))X_a
\\&=&\trace_h \widetilde{R}(d\1\cdot,
\tau(\1))d\1\cdot
\end{eqnarray*}

Thus,
\begin{eqnarray*}
\tau_2(\widetilde{\pi})&=&-\Delta \tau(\widetilde\pi)-
\trace_g\widetilde{R}(d\widetilde{\pi}\cdot,\tau(\widetilde{\pi}))
d\widetilde{\pi}\cdot
\\&=&-\Delta\tau(\1)-\trace_h\widetilde{R} (d\1\cdot,
\tau(\1))di\cdot
\\&=&\tau_2(\1)
\end{eqnarray*}

\end{proof}

\begin{corollary}
Let $\pi:(M,g)\to (N,h)$ be a harmonic Riemannian submersion and
$\widetilde{h}=e^{2\rho}h,\, \rho\in C^\infty(N),$ a conformal
change of metric on $(N,h)$. Denote by
$\widetilde{\pi}=\1\circ\pi:(M,g)\to(N,\widetilde{h})$, where
$\1:(N,h)\to (N,\widetilde{h})$ is the identity map. Then
$\widetilde\pi$ is nonharmonic and biharmonic if and only if $\1$
is nonharmonic and biharmonic.
\end{corollary}

\section{The biharmonicity of the identity map $\1:(N,h)\to (N,\widetilde{h})$}

The previous result allows us to study only the conformal changes
with the property that the identity $\1:(N,h)\to
(N,\widetilde{h})$ is a nonharmonic biharmonic map.

\begin{remark} We note that
$$
\tau(\1)=-e^{2\rho} \tau({\widetilde{\1}})=(2-n)\grad\rho,
$$
so $\1:(N,h)\to(N,\widetilde{h})$ is harmonic if and only if
$\widetilde{\1}:(N,\widetilde{h})\to(N,h)$ is harmonic, and this
holds when $\rho=$ constant for $n>2$. This equivalence is no
longer valid in the case of biharmonicity, (see Example 4.1.(a)),
and the study of the biharmonic equation $\tau_2(\1)=0$ is more
interesting.
\end{remark}

From now on we shall assume $n>2$.

\begin{theorem}
Let $\1:(N,h)\to(N,\widetilde{h})$ be the identity map and
$\widetilde{h}=e^{2\rho}h$ a metric conformally equivalent to $h$.
Then $\1:(N,h)\to(N,\widetilde{h})$ is biharmonic if and only if
\begin{eqnarray}\label{eq: bigrad}
0&=&\trace_h\nabla^2\grad\rho+(2\Delta\rho+(2-n)|\grad\rho|_h^2)\grad\rho
\\&&+\frac{6-n}{2}\grad(|\grad\rho|_h^2)
+\ricci(\grad\rho) \nonumber
\end{eqnarray}
\end{theorem}

\begin{proof}
We recall some well known results about the conformal changes of
metric on a Riemannian manifold. Denote by $\nabla$ the
Levi-Civita connection of the metric $h$ and by
$\widetilde{\nabla}$ the Levi-Civita connection of
$\widetilde{h}.$ The relation between the two connections is given
by $\widetilde{\nabla}_X Y=\nabla_X Y+P(X,Y),$ where the
$(1,2)$-type tensor field $P$ has the expression
\begin{equation}\label{eq:exP}
P(X,Y)=X(\rho)Y+Y(\rho)X-h(X,Y)\grad\rho.
\end{equation}
\\The curvature tensors of $\widetilde\nabla$ and $\nabla$ satisfy
\begin{eqnarray}\label{eq:tenscurb}
\qquad \widetilde{R}^N(X,Y)Z&=&R^N(X,Y)Z+ (\nabla_X
P)(Y,Z)-(\nabla_Y P)(X,Z)
\\&&+P(X,P(Y,Z))-P(Y,P(X,Z)).\nonumber
\end{eqnarray}
Let $\big\{X_a\big\}_{a=1}^n$ be a geodesic frame in $q\in N$. As
$\tau(\1)=(2-n)\grad\rho$ we get
\begin{eqnarray*}
-\Delta\tau(\1)&=&\trace_h\nabla^{\1}\nabla^{\1}\tau(\1)
\\&=&(2-n)\sum_{a=1}^n\big(\widetilde{\nabla}_{X_a}
\widetilde{\nabla}_{X_a}\grad\rho-
\widetilde{\nabla}_{\nabla_{X_a}X_a}\grad\rho\big)
\\&=&(2-n)\sum_{a=1}^n\widetilde{\nabla}_{X_a}
\widetilde{\nabla}_{X_a}\grad\rho.
\end{eqnarray*}
But
\begin{eqnarray*}
\widetilde{\nabla}_{X_a}\grad\rho&=&\nabla_{X_a}\grad\rho
+P(X_a,\grad\rho)
\\&=&\nabla_{X_a}\grad\rho+|\grad\rho|^2X_a
\end{eqnarray*}
and, using \eqref{eq:exP}, we obtain
\begin{eqnarray*}
\widetilde{\nabla}_{X_a} \widetilde{\nabla}_{X_a}
\grad\rho&=&\nabla_{X_a}\nabla_{X_a}\grad\rho+
\frac{3}{2}X_a(|\grad\rho|^2)X_a
+\nabla_{X_a(\rho)X_a}\grad\rho\\&&
-h(X_a,\nabla_{X_a}\grad\rho)\grad\rho.
\end{eqnarray*}
Thus,
\begin{eqnarray*}
-\Delta\tau(\1)&=&(2-n)\big\{\trace_h\nabla^2\grad\rho
+\frac{3}{2}\grad(|\grad\rho|^2)
\\&&+\nabla_{\grad\rho}\grad\rho+(\Delta\rho
+(2-n)|\grad\rho|^2)\grad\rho\big\}.
\end{eqnarray*}

In order to determine $\trace_h
\widetilde{R}(d\1\cdot,\tau(1))d\1\cdot$, by using
\eqref{eq:tenscurb}, we first have
$$
\sum_{a=1}^n P(X_a,P(\grad\rho,X_a))
=|\grad\rho|^2\sum_{a=1}^nP(X_a,X_a)=(2-n)|\grad\rho|^2\grad\rho,
$$
and
$$
\sum_{a=1}^n P(\grad\rho,P(X_a,X_a))
=(2-n)P(\grad\rho,\grad\rho)=(2-n)|\grad\rho|^2\grad\rho.
$$
Also, as
$$
P=d\rho\otimes I+I\otimes d\rho-h\otimes \grad\rho,
$$
where $I:C(TN)\to C(TN)$, $I(X)=X$, $\forall X\in C(TN)$,
\begin{eqnarray*}
(\nabla_XP)(Y,Z)-(\nabla_YP)(X,Z)&=&-h(Z,\nabla_Y\grad\rho)X+h(Z,\nabla_X\grad\rho)Y\\
&&-h(Y,Z)\nabla_X\grad\rho+h(X,Z)\nabla_Y\grad\rho.
\end{eqnarray*}
We shall use the last relation for $X=Z=X_a$ and $Y=\grad\rho$ and
sum. We get
\begin{eqnarray*}
\sum_{a=1}^n\{(\nabla_{X_a}P)(\grad\rho,X_a)-(\nabla_{\grad\rho}P)(X_a,X_a)\}=
-(\Delta\rho)\grad\rho+(n-2)\nabla_{\grad\rho}\grad\rho.
\end{eqnarray*}
Thus,
\begin{eqnarray*}
\trace_h \widetilde{R}(d\1\cdot,\tau(1))d\1\cdot&
=&\trace_h\widetilde{R} (\cdot,\tau(\1))\cdot
\\&=&(2-n)\big\{-\ricci(\grad\rho)
-(\Delta\rho)\grad\rho
\\&&-(2-n)\nabla_{\grad\rho}\grad\rho\big\}.
\end{eqnarray*}
And now, as
\begin{eqnarray*}
\nabla_{\grad\rho}\grad\rho&=&\sum_{a=1}^n
\langle\nabla_{\grad\rho}\grad\rho,X_a\rangle X_a
\\&=&\sum_{a=1}^n\big\{X_a(|\grad\rho|^2)-\langle\nabla_{X_a}
\grad\rho,\grad\rho\rangle\big\}X_a
\\&=&\frac{1}{2}\grad(|\grad\rho|^2).
\end{eqnarray*}
to conclude, by summing, we obtain \eqref{eq: bigrad}.
\end{proof}

\subsection*{} In order to study the connection between the
biharmonicity of $\1:(N,h)\to(N,\widetilde{h})$ and that of
$\widetilde{\1}:(N,\widetilde{h})\to(N,h)$, by using \eqref{eq:
bigrad} and the formula obtained in \cite{PBDK} for
$\tau_2(\widetilde{\1})$, i.e.
\begin{eqnarray*}
\tau_2(\widetilde{\1})&=&(n-2)e^{-4\rho}\big\{\trace_h\nabla^2\grad\rho
+(n-6)\nabla_{\grad\rho}\grad\rho\nonumber\\
&&+2(\Delta\rho-(n-4)\vert
d\rho\vert_h^2)\grad\rho+\ricci(\grad\rho) \big\},
\label{Tau2tild}
\end{eqnarray*}
 we get
$$
\tau_2(\widetilde{\1})+e^{-4\rho}\tau_2(\1)
=e^{-4\rho}(n-2)(n-6)\big\{\grad(\vert\grad\rho\vert^2)
-\vert\grad\rho\vert^2\grad\rho \big\}.
$$
From here,
\begin{proposition} If $n=6$ or
$\grad(\vert\grad\rho\vert^2)=\vert\grad\rho\vert^2\grad\rho$,
then $\widetilde{h}$ is biharmonic with respect to $h$ if and only
if $h$ is biharmonic with respect to $\widetilde{h}$, i.e.
$\tau_2(\1)=0$ if and only if $\tau_2(\widetilde{\1})=0$.
\end{proposition}

\begin{remark}
In general, condition
$$\grad(\vert\grad\rho\vert^2)=\vert\grad\rho\vert^2\grad\rho$$
is not satisfied, (see Example 4.1.(a)).
\end{remark}

In the following we will use the musical isomorphisms
$^\sharp:T^*_qN\to T_qN$, \\ $^\flat:T_qN \to T_q^*N \ (q\in N)$
in order to express condition \eqref{eq: bigrad}.

\begin{proposition}
Let $\1:(N,h)\to(N,\widetilde{h})$ be the identity map and
$\widetilde{h}=e^{2\rho}h$ conformally equivalent to $h$. Then
$\1:(N,h)\to(N,\widetilde{h})$ is biharmonic if and only if
\begin{equation}\label{eq: bidif}
-\Delta\alpha+\frac{6-n}{2}d|\alpha|^2
+(2d^*\alpha+(2-n)|\alpha|^2)\alpha+2(\ricci\
\alpha^\sharp)^\flat=0,
\end{equation}
where $d\rho=\alpha$ (or $\grad\rho=\alpha^\sharp$).
\end{proposition}

\begin{proof}
By using the musical isomorphisms we get \eqref{eq: bigrad}
equivalent to
\begin{eqnarray*}
0=\trace_h\nabla^2\alpha+(2\Delta\rho+(2-n)|\alpha|_h^2)\alpha
+\frac{6-n}{2}d(|\alpha|_h^2) +(\ricci \alpha^\sharp)^\flat.
\end{eqnarray*}

From $\Delta\alpha=(dd^*+d^*d)\alpha=dd^*d\rho$ we obtain
$$(\Delta\alpha)(X)=-(\trace\nabla^2\alpha)(X)+\alpha(\ricci(X))
=(-\trace\nabla^2\alpha+(\ricci\alpha^\sharp)^\flat)(X)$$ and thus
\eqref{eq: bigrad} equivalent to \eqref{eq: bidif}.
\end{proof}

\begin{corollary}
Let $\1:(N,h)\to(N,\widetilde{h})$ be the identity map,
$\widetilde{h}=e^{2\rho}h$ conformally equivalent to $h$ and
denote $d\rho=\alpha$. If $\1$ is biharmonic, then
\begin{equation}\label{eq:consecdif}
(4-n)d|\alpha|^2\wedge\alpha+2(\ricci \
\alpha^\sharp)^\flat\wedge\alpha+d(\ricci \
\alpha^\sharp)^\flat=0.
\end{equation}
\end{corollary}

\begin{proof}
As $\alpha=d\rho$, we have $\Delta\alpha=d(d^*\alpha)$ and
$d\Delta\alpha=0$. By differentiating now \eqref{eq: bidif} it
results
\begin{eqnarray*}
0&=&d\Delta\alpha=d(2d^*\alpha+(2-n)|\alpha|^2)\wedge\alpha+2d(\ricci\
\alpha^\sharp)^\flat
\\&=&2\Delta\alpha\wedge\alpha+(2-n)d|\alpha|^2\wedge\alpha+2d(\ricci\
\alpha^\sharp)^\flat
\\&=&2(4-n)d|\alpha|^2\wedge\alpha+4(\ricci \
\alpha^\sharp)^\flat\wedge\alpha+2d(\ricci \ \alpha^\sharp)^\flat.
\end{eqnarray*}
\end{proof}

\subsection{Biharmonic maps and isoparametric functions}

Equation \eqref{eq: bidif} suggests the study of the relation
between conformal changes which render the identity biharmonic on
Einstein spaces and isoparametic functions.

Let $(N,h)$ be a Riemannian manifold and $f\in C^{\infty}(N)$.
\begin{definition}\cite{PB}.
$f:N\to \mathbb{R}$ is called {\it isoparametric} if there
are smooth real functions $\gamma$ and $\sigma$ such that

\begin{equation}
\left\{
\begin{array}{l}
(i)\,\ |df|^2=\gamma\circ f,
\\ \mbox{} \\
(ii)\,\ \Delta f=\sigma\circ f,
\end{array}
\right.\label{eq:izoparam1}
\end{equation}
at all points where $\grad f\neq 0.$
\end{definition}

\begin{remark} $f$ is isoparametric if and only if $\forall
q\in N$ with $\grad_q f\neq 0,$
\begin{equation}
\left\{
\begin{array}{l}
(i)\,\ \grad(|\grad f|) \textrm{ is parallel to} \grad f,
\\ \mbox{} \\
(ii)\,\ \grad(\Delta f) \textrm{ is parallel to} \grad f.
\end{array}
\right.\label{eq:izoparam2}
\end{equation}

Indeed, from \eqref{eq:izoparam1}(i), by using $|df|^2=|\grad
f|^2$, we get
\begin{align*}
\grad(|df|^2)=&2|df|\grad(|\grad f|)\\
=&(\gamma'\circ f)\grad f.
\end{align*}
Thus, in any point with $\grad f\neq 0$ we have
$$
\grad(|\grad f|)=\frac{\gamma'\circ f}{2|df|}\grad f,
$$
so $\grad(|\grad f|)$ is parallel to $\grad f$. Also,
\eqref{eq:izoparam1}(ii) implies
$$\grad(\Delta f)=(\sigma'\circ f)\grad f.$$

Conversely, let $g$ be a smooth function on $N$. Condition $\grad
f$ parallel to $\grad g$ away from its zeros insures the existence
of $ \overline{g}:f(M)\to g(M)$, with $\overline{g}\circ f=g$.
This remark, together with \eqref{eq:izoparam2}(i) leads, for
$g=|\grad f|$, to \eqref{eq:izoparam1}(i) and together with
\eqref{eq:izoparam2}(ii), for $g=\Delta f$, to
\eqref{eq:izoparam1}(ii).

\end{remark}

\begin{remark}\cite{PB}.
Condition \eqref{eq:izoparam1}(i) (or \eqref{eq:izoparam2}(i)) is
equivalent to the fact that the regular levels of $f$ are
parallel, i.e. integral curves of $\grad f/|\grad f|$ are
geodesics, and condition \eqref{eq:izoparam1}(ii) (or
\eqref{eq:izoparam2}(ii)) with the fact that the regular levels of
$f$ have constant mean curvature.
\end{remark}

\begin{lemma}
Let (N,h) be an Einstein space. If $\rho$ satisfies
\eqref{eq:izoparam2}(i) then \eqref{eq: bidif} implies
$$\Delta\alpha-F(q,\alpha)\,\alpha=0,$$
with $F$ a real valued function depending on $q$ and $\alpha$.
\end{lemma}
\begin{proof}
As $N$ is an Einstein space we have
$$(\ricci \ \alpha^\sharp)^\flat=c \,\alpha,$$
where $c$ is the Einstein constant of $N$, and (i) implies
$$d|\alpha|^2=d\big(\gamma\circ \rho)=(\gamma'\circ \rho\big)\alpha,$$
so
\begin{eqnarray*}
\Delta\alpha&=&(2d^*\alpha+(2-n)|\alpha|^2
+\frac{6-n}{2}(\gamma'\circ\rho)+c)\,\alpha
\\&=&F(q,\alpha)\,\alpha.
\end{eqnarray*}

\end{proof}

\begin{corollary}
Let (N,h) be an Einstein space and $\rho:N\to \mathbb{R}$
satisfying \eqref{eq:izoparam2}(i). If $\rho$ satisfies \eqref{eq:
bidif}, then $\rho$ is isoparametric.
\end{corollary}
\begin{proof}
We will show that $\Delta\rho$ is constant along the  integral
curves of any vector field tangent to the level surfaces of
$\rho$, condition equivalent to \eqref{eq:izoparam2}(ii). Indeed,
\begin{eqnarray*}
X(\Delta\rho)&=&X(d^*d\rho)=d(d^*d\rho)(X)
\\&=&\Delta\alpha(X)=F(q,\alpha)\alpha(X)=F(q,\alpha)X(\rho)\\
&=&0,
\end{eqnarray*}
$\quad\forall X$ tangent to the regular levels of $\rho$.
\end{proof}

\begin{theorem}
Let $(N^n,h), n\neq 4$, be an Einstein space and
$\widetilde{h}=e^{2\rho}h$ a metric conformally equivalent to $h$
such that $\1:(N,h)\to (N,\widetilde{h})$ is biharmonic. Then
$\rho:N\to\mathbb{R}$ is isoparametric.
\\Conversely, if $f:N\to\mathbb{R}$ is isoparametric, then away
from critical points of $f$, there exists a reparametrisation
$\rho=\rho\circ f$ such that $\1:(N,h)\to (N,\widetilde{h})$ is
biharmonic.
\end{theorem}

\begin{proof}
We will prove that, in these conditions, $\rho$ satisfies
\eqref{eq:izoparam2}(i) and so, the conclusion follows from the
last corollary. Equation \eqref{eq:consecdif} reduces to
$$(n-4)d(|d\alpha|^2)\wedge\alpha=0.$$
As $n\neq 4$, this implies $X(|\grad\rho|^2)=0,\,\forall X
\,\textrm{with}\, X(\rho)=0,$ and thus the part of $\grad(|\grad
\rho |^2)$ tangent to the regular levels of $\rho$ is zero and
\eqref{eq:izoparam2}(i) follows.

Conversely, if $f:N\to\mathbb{R}$ is a isoparametric non-constant
function, we will determine a reparametrisation $\rho=\rho\circ f$
of $f$ with $\alpha=d\rho$ satisfying \eqref{eq: bidif}.

If $|df|^2=\gamma\circ f$, by choosing $s$ solution of the
differential equation $s'=\frac{1}{\gamma^{\frac{1}{2}}}$, we can
reparametrize $f$ such that $s=s\circ f$, and $|ds|^2=1$. The
function $s$ is isoparametric with $\Delta s=\sigma\circ s$, for a
smooth function $\sigma$. We have
$$d\rho=(\rho'\circ s)ds,$$
$$\Delta\rho=(-\rho''+\rho'\sigma)\circ s.$$
From here we get
$$\Delta\alpha=dd^*d\rho=d(\Delta\rho)=
(-\rho'''+\rho''\sigma+\rho'\sigma')ds,$$
$$d|\alpha|^2=2\rho'\rho''ds,$$
$$(\ricci \ \alpha^\sharp)^\flat=c\,\alpha=c\rho'ds.$$
and, by substituting into \eqref{eq: bidif}, we obtain a second
order ordinary differential equation in $y(s)=\rho'(s)$
\begin{equation}\label{eq:dif}
y''-\sigma y'+(4-n)yy'+(2c-\sigma')y+2\sigma y^2+(2-n)y^3=0.
\end{equation}
The general theory of differential equations offers local
solutions for this type of problem. Finally, $\rho(s)=\int
y(s)ds$.
\end{proof}

\section{Examples}

In this section we present two examples of isoparametric functions
and their reparametrisation such that the identity map
$\1:(N,h)\to (N,\widetilde{h})$ is biharmonic.
\begin{example}
Let $N=\mathbb{R}^n_+=\{x\in \r^n/ x^i>0, i=\overline{1,n}\}$
endowed with the canonic metric and $s:N\to\mathbb{R},\, s(x^1,
x^2,...,x^n)=\sum_{i=1}^n \alpha^ix^i,$ where
$\sum_{i=1}^n(\alpha^i)^2=1, \alpha^i>0$. Note that $s$ is
isoparametric with $|ds|^2=|\alpha|^2=1$ and $\Delta s=0 $.
Equation \eqref{eq:dif} becomes
\begin{equation}\label{eq: ex1}
y''+(4-n)yy'+(2-n)y^3=0.
\end{equation}

We will search special solutions satisfying $y'=ay^2$. We have
$y''=2ayy'$, and \eqref{eq: ex1} is equivalent, excluding $y=0$,
with
$$
2a^2+(4-n)a+(2-n)=0.
$$
We get $a=-1$ and $a=\frac{n-2}{2}$.

a) If $a=-1$, then $\rho(s)=\ln s$ and $\widetilde{h}=s^{2}h$.

Also, for $s=x^1$ we have
$$
\grad(|\grad\rho|^2)=\grad\frac{1}{(x^1)^2}=-\frac{2}{(x^1)^3}e_1,
$$
and
$$
|\grad\rho|^2\grad\rho=\frac{1}{(x^1)^3}e_1,
$$
so $h$ is biharmonic respect to $\widetilde{h}$, but not
conversely.

b) If $a=\frac{n-2}{2}$, then $\widetilde{h}=s^{\frac{4}{2-n}}h$.
\end{example}

\begin{example}
Let $N=\mathbb{R}^n$ be endowed with the canonic metric and
$s:N\to\mathbb{R},\, s(x^1, x^2,...,x^n)
=\sqrt{(x^1)^2+...+(x^n)^2}$. $s$ is isoparametric with $|ds|^2=1$
and $\Delta s=-(n-1)/s $, so $\sigma(s)=-(n-1)/s$. In this case
\eqref{eq:dif} becomes
\begin{equation}
y''+\frac{n-1}{s} y'+(4-n)yy'-\frac{n-1}{s^2}y-2\frac{n-1}{s}
y^2+(2-n)y^3=0. \label{eq: ex2}
\end{equation}
For this, we will search for solutions satisfying $y=a/s$, with
$a\in \r$, constant. \eqref{eq: ex2} becomes, excluding $y=0$,
$$
(2-n)a^2-(n+2)a+4-2n=0.
$$
The discriminant of this equation is positive only for $n \in\{1,
2, 3, 4\}$, and $n\neq 2, 4$ implies:

a) for $n=1$, $a\in \{1,2\}$;

b) for $n=3$, $a\in\big\{\frac{5\pm\sqrt{17}}{2}\big\}$.
\end{example}

\end{document}